\documentclass[12pt]{pmsart}
\allowdisplaybreaks
\newcommand{\NN}{ {\bf N} }

\newcommand{\CC}{{\mathbb C}}

\newcommand{\ff}{\varphi}
\newcommand{\la}{\langle}
\newcommand{\ra}{\rangle}
\newcommand{\kk}{{\kappa}}
\newcommand{\Tr}{\text{Tr}}
\newcommand{\tr}{\text{tr}}
\newcommand{\id}{\text{id}}

\newcommand{\cA}{\mathcal{A}}

\newcommand{\ab}{\allowbreak}
\newcommand{\ds}{\displaystyle}

\ShortAuthors{Mingo \& Speicher}
\ShortTitle{Schwinger-Dyson Equations}
\firstpage{1}
\volume{XX}
\fasc{2}
\years{2013}
\pages{1-13}
\dedicated{}
\received{\today}
\lastrevision{\today}

\begin{document}

\title{Schwinger-Dyson Equations: Classical and Quantum}
\NumberOfAuthors{2}

\FirstAuthor{James A. Mingo}

\SecondAuthor{Roland Speicher}

\FirstAuthorAffiliation{Queen's University}
\FirstAuthorAddress{Department of Mathematics and Statistics\\
                    Kingston, Ontario\\ K7L 3N6, Canada}
\FirstAuthorCity{Kingston}
\FirstAuthorURL{www.mast.queensu.ca/{\raise.17ex\hbox{%
$\scriptstyle\sim$}}mingo}
\FirstAuthorEmail{mingo@mast.queensu.ca}
\FirstAuthorThanks{Supported by the Natural Sciences and
  Engineering Research Council of Canada}

\SecondAuthorAffiliation{Universit\"at des Saarlandes}
\SecondAuthorAddress{FR $6.1-$Mathematik\\ Postfach 15 11 50\\ 
D-66123 Saarbr\"{u}cken, Germany}
\SecondAuthorCity{Saarbr\"ucken}
\SecondAuthorURL{math.uni-sb.de/ag/speicher}
\SecondAuthorEmail{speicher@math.uni-sb.de}
\SecondAuthorThanks{Supported by the Alfried Krupp von
  Bohlen und Halbach Stiftung (R\"uckkehr deutscher
  Wissenschaftler aus dem Ausland) and by the DFG (SP
  419/8-1)}

\keywords{free probability, random matrices, Schwinger-Dyson equation}
\MSCcodes{60B20, 46L54}
\maketitle

\section{Introduction}
In this note we want to have another look on Schwinger-Dyson
equations for the eigenvalue distributions and the
fluctuations of classical unitarily invariant random matrix
models. We are exclusively dealing with one-matrix models,
for which the situation is quite well understood. Our point
is not to add any new results to this, but to have a more
algebraic point of view on these results and to understand
from this perspective the universality results \cite{AJM,J}
for fluctuations of these random matrices.  We will also
consider corresponding non-commutative or ``quantum" random
matrix models and contrast the results for fluctuations and
Schwinger-Dyson equations in the quantum case with the
findings from the classical case.

\section{Notations and Prerequisites}
\subsection{Free probability theory}
For the basic notions and results about free probability
theory we refer to the books \cite{VDN,NS}; in particular,
we will follow the latter in regard of the definitions and
fundamental results on free cumulants.

\subsection{Non-commutative Derivatives}
We will denote by $\partial$ and $D$ the non-commutative and
the cyclic derivative, respectively; see, for example,
\cite{V-cyclic} for definitions and basic properties; note
that in \cite{V-cyclic} the cyclic derivative is denoted by
$\delta$. We will only use these derivatives in the
one-variable case; then, the cyclic derivative $D$ coincides
with usual differentiation.  On the algebra $\CC\la x\ra$ of
polynomials in one variable $x$ these derivatives are given
by
\begin{align*}
D: \CC\la x\ra &\to \CC\la x\ra\\
x^n&\mapsto Dx^n:=nx^{n-1}
\end{align*}
and
\begin{align*}
\partial:\CC\la x\ra &\to\CC\la x\ra\otimes \CC\la x\ra\\
x^n&\mapsto \partial x^n:=\sum_{k=0}^{n-1} x^k\otimes x^{n-k-1}
\end{align*}

\subsection{The Chebyshev polynomials}

We will use the Chebyshev polynomials of first and second
kind, for the interval $[-2, 2]$. The ones orthogonal with
respect to the semicircle (second kind) are denoted by
$S_n$, the ones orthogonal with respect to the arc-sine
distribution (first kind) by $C_n$; compare \cite{KMS}.  We
have
$$C_0(x)=2,\qquad C_1(x)=x,\qquad C_2(x)=x^2-2,\qquad C_3(x)=x^3-3x$$
and
$$xC_n(x)=C_{n+1}(x)+C_{n-1}(x)\qquad (n>1);$$
and
$$S_0(x)=1,\qquad S_1(x)=x,\qquad S_2(x)=x^2-1,\qquad
S_3(x)=x^3-2x$$
and
$$xS_n(x)=S_{n+1}(x)+S_{n-1}(x)\qquad (n>1).$$
One has, for $n\geq 0$, the the following identities:
$$D C_n=n S_{n-1},\qquad
\partial S_n=\sum_{k=0}^{n-1}S_k\otimes S_{n-k-1}$$
Furthermore, $C_n=S_n-S_{n-2}$ (those are true for all $n\geq 0$, if we set 
$S_{-2}(x)\ab =  -1$ and $S_{-1}(x) = 0$) and for $n,m\geq 0$
\begin{align*}
S_nS_m&=S_{n+m}+S_{n+m-2}+\cdots+S_{\vert n-m\vert}\\
C_nC_m&=C_{n+m}+C_{\vert n-m\vert}
\end{align*}

These imply that we have for all $n,m\geq 0$
\begin{equation}\label{equation:first-times-second}
C_nS_m=\begin{cases}
S_{n+m}+S_{m-n},& n\leq m\\
S_{n+m},& n=m+1\\
S_{n+m}-S_{n-m-2},& n\geq m+2
\end{cases}.\end{equation}

\subsection{Non-commutative probability space of second order}
A \emph{second order non-commutative probability space}
$(\cA,\ff_1,\ff_2)$ consists of a unital algebra $\cA$, a
tracial linear functional $\ff_1:\cA\to\CC$ with $\ff(1)=1$
and a bilinear functional $\ff_2:\cA\times\cA\to\CC$, which
is symmetric in both arguments, i.e.,
$\ff_2(a,b)=\ff_2(b,a)$ for all $a,b\in\cA$, tracial in each
of its both arguments and which satisfies
$\ff_2(a,1)=0=\ff_2(1,b)$ for all $a,b\in\cA$. Compare
\cite{MSS} for more information.

\section{Schwinger-Dyson equations for classical
            unitarily invariant ensembles}

We will be interested in unitarily invariant random
matrices; the most prominent class of random matrices of
this type is given by a density of the following form.  We
consider Hermitian $N\times N$-random matrices
$A=(a_{ij})_{i,j=1}^N$ equipped with the probability measure
\begin{equation}\label{ensemble}
d\mu_N(A)=\frac 1{Z_N}\exp\bigl\{-N\Tr[P(A)]\bigr\}dA,
\end{equation}
where
$$dA=\prod_{1\leq i<j\leq N}d\,\text{Re}\, a_{ij}\,
d\,\text{Im}\, a_{ij} \prod_{i=1}^N da_{ii}.$$ 
Here, $P$ is a polynomial in one variable, which we will
address in the following as ``potential", and $Z_N$ is a
normalization constant to make \eqref{ensemble} into a
probability distribution.

At least formally, it is quite easy to see that the
asymptotic eigenvalue distribution and fluctuations of these
ensembles satisfy in the large $N$-limit the following
so-called Schwinger-Dyson equations (see \cite[Chapter
  8]{G1}, also called the method of equation of motion or
the loop equation in \cite[Chapter 6]{E}). .  We will ignore
all analytic questions and just work in the algebraic
setting; thus we take our non-commutative probability space
$\cA=\CC\la x\ra$ as the polynomials in one variable $x$.

\begin{definition}
Let $(\CC\la x\ra,\ff_1,\ff_2)$ be a non-commutative
probability space of second order and $V\in\CC\la x\ra$ a
polynomial in $x$.  We put $\xi:=D V(x)\in\CC\la x\ra$.  We
say that $\ff_1$ satisfies the \emph{first order
  Schwinger-Dyson equations} for the potential $V$ if we
have for all $p(x)\in\CC\la x\ra$
\begin{equation}\label{SW:first}
\ff_1\bigl(\xi p(x)\bigr)=\ff_1\otimes \ff_1\bigl(\partial p(x)\bigr)
\end{equation}
(i.e., $\xi$ is the conjugate variable for $x$).
If we have in addition that for all $p(x),q(x)\in\CC\la x\ra$
\begin{eqnarray}\label{SD:second}\lefteqn{
\ff_2(\xi p(x), q(x))} \\
& = &
\ff_2\bigl([\ff_1\otimes\id+\id\otimes \ff_1](\partial p(x)), q(x)\bigr)+\ff_1(p(x)D q(x)), \notag
\end{eqnarray}
then $(\ff_1,\ff_2)$ satisfies the \emph{second order
  Schwinger-Dyson equations}.
\end{definition}

Corresponding analogues exist also for the case of several
matrices, but since we have nothing substantial to say about
the multi-variate case we will stick in the following to the
one-matrix case. Existence and uniqueness of the solution of
these equations (under positivity requirements for $\ff_1$)
are well-studied in the one-matrix case, and are one of the
main problems in random matrix theory for the case of
several variables; for some positive results in the latter
case see \cite{Gui}.

We will in the following ignore the uniqueness question and
present a solution to the Schwinger-Dyson equations for the
one-matrix case.

\begin{theorem}
For a given $V\in \CC\la x\ra$, we decompose $D V$ with
respect to the Chebyshev polynomials of the first kind
$$\xi=D V(x)=\sum_{n\geq 0} \alpha_n C_n(x).$$ Assume that
we have normalized $V$ in such a way that $\alpha_0=0$ and
$\alpha_1=1$.  We define on $\CC\la x\ra$ a $\ff_1$ by
$$\ff_1(S_{n}(x)):=\alpha_{n+1} \qquad (n\geq 0)$$
(note that we need $\ff_1(1)=\alpha_1=1$ for this)
and a $\ff_2$ by
$$\ff_2(C_n(x),C_m(x)):=n\delta_{nm}\qquad (n,m\geq 0).$$
Then $\ff_1$ and $\ff_2$ satisfy the first and second order
Schwinger-Dyson equations for the potential $V$.
\end{theorem}

The prescriptions above provide well-defined and unique
$\ff_1$ and $\ff_2$, because both $\{S_n\mid n\geq 0\}$ and
$\{C_n\mid n\geq 0\}$ are linear bases of $\CC\la x\ra$.

Note also the crucial fact that $\ff_2$ does not depend on
$V$. Actually, our definition of $\ff_2$ is in essence just
a reformulation of the universality of the asymptotic
fluctuations for the random matrix ensemble given by
\eqref{ensemble}.  In the physical literature this
observation goes at least back to Politzer \cite{P},
culminating in the paper of Ambj\o{}rn et al. \cite{AJM},
whereas a proof on the mathematical level of rigour is due
to Johansson \cite{J}. The above theorem arouse out of our
attempts to understand this universality result. Actually,
it can (and should) also be seen as a streamlined algebraic
proof of this universality result.

Our original motivation in this context was to look for
multivariate versions of this result. As will be seen from
the following proof, the result relies crucially on various
algebraic properties of the Chebyshev polynomials, for which
no multivariate version exists. Thus it should be clear that
the universality result is a genuine one-dimensional
phenomena. Actually, in \cite{MSS} we have shown, by using
the machinery of second order freeness, that for one of the
most canonical families of several random matrices the
fluctuations depend indeed on the potential $V$.

\begin{proof}
Consider the first order. We have to show that
$$\ff_1(\xi p(x))=\ff_1\otimes \ff_1(\partial p(x))$$ 
for all $p(x)\in\CC\la x\ra$. By linearity, it suffices to
treat the cases $p(x)=S_m(x)$ for all $m\geq 0$. So fix such
an $m$.  Thus we have to show
$$\sum_{n\geq 0} \alpha_n
\ff_1(C_n(x)S_m(x))=\ff_1\otimes\ff_1(\partial S_m(x))$$
For the left hand side 
we have
\begin{align*}
\sum_n\alpha_n \ff_1(C_nS_m)
&=\sum_{n\leq m}
\alpha_{n}\bigl(\ff_1(S_{n+m})+\ff_1(S_{m-n})\bigr)+\alpha_{m+1}\ff(S_{2m+1})\\
&\qquad+\sum_{n\geq
 m+2}\alpha_n\bigl(\ff_1(S_{n+m})-\ff_1(S_{n-m-2})\bigr)\\
&=\sum_{n\leq m}
\alpha_{n}\bigl(\alpha_{n+m+1}+\alpha_{m-n+1}\bigr)+\alpha_{m+1}\alpha_{2m+2}\\
&\qquad+\sum_{n\geq
 m+2}\alpha_n\bigl(\alpha_{n+m+1}-\alpha_{n-m-1}\bigr)\\
 &=\sum_{n}
\alpha_{n}\alpha_{n+m+1}-\sum_{n\geq
 m+2}\alpha_n\alpha_{n-m-1}+
 \sum_{n\leq m}
\alpha_{n}\alpha_{m-n+1}.
\end{align*}
But the first two sums cancel as the summation in $n$ starts
at $n=1$ (because $\alpha_0=0$), and thus we remain with
exactly the same as in
$$\ff_1\otimes\ff_1(\partial S_m(x))=\sum_{k=0}^{m-1}\ff_1(S_k)\ff_1(S_{m-k-1})=
\sum_{k=0}^{m-1} \alpha_{k+1}\alpha_{m-k}.$$

Now consider the second order. For this we have to show that
\begin{eqnarray*}\lefteqn{
\ff_2(\xi p(x), q(x)) } \\
& = &
\ff_2([\ff_1\otimes\id+\id\otimes \ff_1](\partial p(x)), q(x))+\ff_1\bigl(p(x) \cdot D q(x)\bigr)
\end{eqnarray*}
for all $p$ and $q$. Again, by linearity, it is enough to
show this for $p=C_m$ and $q=C_k$, for arbitrary $m,k\geq
0$. Thus we have to show
\begin{eqnarray}\label{equation:modified-second-order-schwinger-dyson}\lefteqn{
\sum_{n\geq 0}\alpha_n \ff_2(C_n C_m,C_k)} \\
& = &
\ff_2([\ff_1\otimes\id+\id\otimes \ff_1](\partial C_m), C_k)+\ff_1(C_m k S_{k-1}).\notag
\end{eqnarray}
We have (note that we set $S_{-2}=-1$ and $S_{-1}=0$)
\begin{align*}
\partial C_m&=\partial (S_m-S_{m-2})\\&=\sum_{l=0}^{m-1} S_l\otimes S_{m-l-1} -\sum_{l=0}^{m-3} S_l\otimes 
S_{m-2-l-1}\\
&=\sum_{l=0}^{m-1} S_l\otimes \tilde C_{m-l-1},
\end{align*}
where $\tilde C_r=C_r$ for $r\geq 1$ and $\tilde C_0=1=S_0$.
Thus we have
$$\ff_1\otimes \id (\partial C_m)=\sum_{l=0}^{m-1} \ff_1(S_l) \tilde C_{m-l-1}=\sum_{l=0}^{m-1} \alpha_{l+1} \tilde C_{m-l-1}.$$
Hence
\begin{align}\label{equation:RHS(ii)}
\ff_2([\ff_1\otimes\id +\id\otimes \ff_1](\partial C_m), C_k)\qquad\qquad&\\
=2\sum_{l=0}^{m-1} \alpha_{l+1}
\ff_2(\tilde C_{m-l-1},C_k)=
\begin{cases}
2 \alpha_{m-k} k, & k\leq m\\
0,& k>m
\end{cases}.\notag
\end{align}
Next using the formula (\ref{equation:first-times-second})
for $C_m S_{k-1}$ we have
\[
\ff_1(C_m k S_{k-1})
=
\begin{cases}
\alpha_{m+k} + \alpha_{k-m} & m \leq k - 1\\
\alpha_{m+k}                & m = k \\
\alpha_{m+k} - \alpha_{m-k} & m  > k
\end{cases}
\]
If we add this to the right hand side of
(\ref{equation:RHS(ii)}) we get that the right hand side of
(\ref{equation:modified-second-order-schwinger-dyson}) is
$k(\alpha_{m+k} + \alpha_{|m-k|})$. Finally let us check the
left hand side of
(\ref{equation:modified-second-order-schwinger-dyson}).
\begin{align*}
\sum_{n\geq 0} \alpha_n \ff_2(C_m C_n, C_k)
& = 
\sum_{n \geq 1} \alpha_n \{ \ff_2(C_{m+n}, C_k) + \ff_2(C_{|m-n|}, C_k)\} \\
& = 
k \begin{cases}
\alpha_{m+k} + \alpha_{k-m} & m < k \\
\alpha_{m+k} + \alpha_{m-k} & m \geq k
\end{cases} \\
& = 
k(\alpha_{m+k} + \alpha_{|m-k|})
\end{align*}
Thus both sides of
(\ref{equation:modified-second-order-schwinger-dyson}) equal
$k(\alpha_{m+k} + \alpha_{|m-k|})$ as claimed.
\end{proof}

\section{Quantum matrix models}

Now we want to consider non-commutative (or ``quantum")
analogues of our classical random matrix models; i.e., we
consider matrices where the entries are not commutative
random variables, but in general non-commu\-tative ones. We
want to address the question about fluctuations in such a
context.

The essential property of the classical ensemble
\eqref{ensemble} is the invariance under unitary
conjugation, i.e., the joint distribution of the entries of
$A=(a_{ij})_{i,j=1}^N$ does not change if we go over to the
conjugated matrix $B:=UAU^*$ for any $N\times N$ unitary
matrix $U$. We will now look on analogues of this for
quantum $N\times N$ matrices $A=(a_{ij})_{i,j=1}^N$ (where
the entries $a_{ij}$ come from some non-commutative
probability space $(\cA, \ff)$), but where we ask not just
for invariance under conjugation by classical unitary
matrices, but -- in line with the idea that one should also
replace classical symmetries by corresponding quantum
symmetries in a non-commutative context -- for the stronger
corresponding invariance under the action of the quantum
unitary group $U_N^+$. By \cite{CS}, a big class of such
invariant matrices are given by the requirement that $A$ is
free from $M_N(\CC)$.  Another characterization of this is
as follows: the matrix $A$ is $R$-cyclic (in the sense of
\cite{NSS}) and the non-vanishing cumulants of its entries
depend only on the length of the cumulant.  A way to
construct such quantum random matrices is by compressing
some random variable $a$ with free matrix units; compare
Lecture~14 in \cite{NS}.

Recall that a matrix $A = (a_{ij})_{i,j=1}^N \in M_N(\cA)$
is $R$-cyclic if for every $n$ we have $\kk_n(a_{i(1)j(1)},
\dots, a_{i(n)j(n)}) = 0$ unless $j(1) = i(2), \dots, j(n) =
i(1)$ (see \cite[Lecture 20]{NS}). Suppose we have a family
of matrices $\{A_1, \dots , A_s\}$, where we write $A_k =
(a^{(k)}_{ij})_{i,j=1}^N$. The family is $R$-cyclic if for
every $n$ and for every $r(1), \dots, r(n)$ we have
$\kk_n(a^{(r(1))}_{i(1)j(1)}, \dots, a^{(r(n))}_{i(n)j(n)})
= 0$ unless $j(1) = i(2), \dots, j(n) = i(1)$. In
\cite[Theorem 4.3]{NSS} it was shown that matrices from the
algebra generated by a $R$-cyclic family are themselves
$R$-cyclic (see also \cite[Exercise 20.23]{NS}).

So let us in the following fix a selfadjoint random variable
$a$ and denote by $\kk_n:=\kk_n(a,\dots,a)$ the free
cumulants of $a$.  Then, for each $N\in\NN$, we define a
quantum random matrix $A=(a_{ij})_{i,j=1}^N$ by prescribing
the free cumulants of the entries as follows: the cyclic
cumulants of the matrix entries are given by
\begin{equation}
\kk_n(a_{i(1)i(2)},\dots,a_{i(n)i(1)})=\frac 1{N^{n-1}}\kk_n(a,\dots,a),
\end{equation}
all other cumulants being zero.

We are interested in calculating, for $N\to\infty$,
cumulants of traces of powers of $A$.  Fix $n\geq 1$ and
$k(1),\dots,k(n)\geq 1$.  Let $k = k(1) + \cdots + k(n)$. We
have
\begin{eqnarray}\label{equation:nthcumulant}\lefteqn{\kern3em
\kk_n(\Tr(A^{k(1)}),\dots,\Tr(A^{k(n)}))}\\
& = & \kern-1.5em
\sum_{i(1), \dots, i(k)=1}^N
\kk_n(a_{i(1)i(2)} \cdots a_{i(k_1)i(1)}, 
a_{i(k_1+1)i(k_1+2)} \cdots a_{i(k_1+k_2)i(k_1+1)}, \dots, \notag\\ 
&& \qquad
a_{i(k_1+\cdots +k_{n-1} +1)i(k_1+\cdots +k_{n-1} +2)} \cdots a_{i(k_1+\cdots +k_n)i(k_1+\cdots +k_{n-1} +1)}) \notag
\end{eqnarray}
Now since $A$ is $R$-cyclic, the family \{$A^{k(1)}, \dots,
A^{k(n)}\}$ is an $R$-cyclic family; so we know that only
cyclic cumulants in these powers are different from
zero. This means that in the sum above only terms with
$i(1)=i(k_1+1)=\dots=i(k_1 + \cdots +k_{n-1}+1)$ can be
different from 0.

Next we use the formula for cumulants with products as
entries (see \cite[Lecture 11]{NS}) and write 

\noindent
$\ds
\kk_n(a_{i(1)i(2)} \cdots a_{i(k_1)i(1)}, 
a_{i(k_1+1)i(k_1+2)} \cdots a_{i(k_1+k_2)i(k_1+1)}, \dots,
\hfill$

$\hfill\ds
\qquad
a_{i(k_1+\cdots +k_{n-1} +1)i(k_1+\cdots +k_{n-1} +2)} \cdots a_{i(k_1+\cdots +k_n)i(k_1+\cdots +k_{n-1} +1)}) 
$

\noindent
as
\begin{multline*}
\sum_{\pi} \kk_\pi(
a_{i(1)i(2)}, \dots a_{i(k_1)i(1)}, \dots,\\
a_{i(k_1+\cdots +k_{n-1} +1)i(k_1+\cdots +k_{n-1} +2)}, \dots,
 a_{i(k_1+\cdots +k_n)i(k_1+\cdots +k_{n-1} +1)})
\end{multline*}
where the sum runs over all $\pi\in NC(k(1)+\dots+k(n))$
which have the property that they connect the blocks
of $$\tau = \{(1, \dots, k(1)), \dots, (k(1) + \cdots +
k(n-1)+1, \dots,\ab k(1) + \cdots + k(n))\}.$$ In the
language of \cite[Definition 9.15]{NS} this means that $\pi
\vee \tau = 1_k$.

For such a $\pi$ to make a non-zero contribution some
relations on the indices must be satisfied. Let us work out
what this means. Recall that there is an embedding of
$NC(k)$ into $S_k$ the symmetric group on $[k]$; namely put
the elements of the blocks of $\pi\in NC(k)$ in increasing
order and regard them as the cycles of permutation (see
e.g. \cite[Remark 23.24]{NS}).

Suppose $(j_1, \dots, j_r)$ is a block of $\pi$, then the
corresponding factor of $\kk_\pi$ is
\[
\kk_r(a_{i(j_1)i(j_1+1)}, \dots, a_{i(j_r)i(j_r+1)}).
\]
In order for this cumulant to be different from 0 we must have
\[
i(j_1+1) = i(j_2), i(j_2 + 1) = i(j_3), \dots, 
i(j_r + 1 ) = i(j_1).
\]
Let $\gamma \in S_k$ be the permutation with the single
cycle $(1, \dots, k)$. Then our relation on $i$ can be
expressed as
\[
i(j_k) = i(j_{k-1} + 1) = i(\gamma(j_{k-1})) = \ab i(\gamma(\pi^{-1}(j_k)))
\] 
or that $i = i\circ \gamma\pi^{-1}$. An important fact of
the embedding of $NC(k)$ into $S_k$ is that the Kreweras
complement of $\pi$, $K(\pi) = \pi^{-1}\gamma$. What we have
here is the `other' Kreweras complement $\gamma\pi^{-1}$
which is the conjugation of $K(\pi)$ by $\gamma$ (see
\cite[Exercise 9.23(1)]{NS}).

\begin{figure}%
{In this example $k(1) = 4$, $k(2) = 2$, $k(3)= 3$, and
  $k(4) = 4$. The partition $\pi = \{
  (1,10,13)(2,4,5,9)(3)(6)(7,8)(11,12) \}$; the `other'
  Kreweras complement of $\pi$ is $\gamma\pi^{-1} =
  \{(\textbf{1})(2,
  \textbf{10})(3,4)(\textbf{5})(6,\textbf{7},9)(8)(11,13)(12)
  \}$. Note that since $\pi \vee \tau = 1_{13}$, each of the
  points of $\{1, 5,7, 10\}$ is in a separate block of
  $\gamma\pi^{-1}$. There are $d = \#(\gamma\pi^{-1}) - n +
  1 = 5$ degrees of freedom in $i(1), \dots, i(13)$, namely
  $i(1) = i(2) = i(5) = i(6) = i(7) = i(9) = i(10)$, $i(3) =
  i(4)$, $i(11) = i(13)$, $i(8)$, and $i(12)$; i.e. we join
  the blocks of $\gamma\pi^{-1}$ containing a bold number
  and the rest remain as they are.}{}
\includegraphics{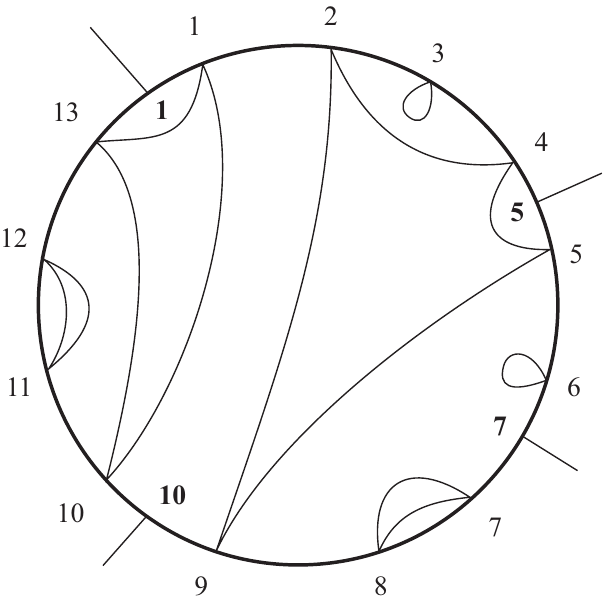}
\end{figure}

Thus in order for
\begin{multline}\label{equation:nonvanishing}
\kk_\pi(
a_{i(1)i(2)}, \dots a_{i(k_1)i(1)}, \dots,
a_{i(k_1+\cdots +k_{n-1} +1)i(k_1+\cdots +k_{n-1} +2)}, \dots,\\
 a_{i(k_1+\cdots +k_n)i(k_1+\cdots +k_{n-1} +1)}) \not= 0
\end{multline}
we must have that $i$ is constant on the cycles of
$\gamma\pi^{-1}$. This is true for any $\pi \in NC(k)$. Let
us now consider what happens when we add the condition $\pi
\vee \tau = 1_k$. According to \cite[Lemma 14]{MST} $\pi
\vee \tau = 1_k$ if and only if each point of the set
$\{k(1), k(1) + k(2), \dots, k(1) + \cdots + k(n)\}$ lies in
a different block of $K(\pi)$; after conjugation by $\gamma$
this condition becomes that each point of $\{1, k(1) + 1,
\dots, k(1) + \cdots + k(n-1) + 1\}$ is in a separate cycle
of $\gamma\pi^{-1}$. Now recall that we had earlier observed
that $R$-cyclicity forced us to have
$i(1)=i(k_1+1)=\dots=i(k_1 + \cdots +k_{n-1}+1)$ in order
for the corresponding term of (\ref{equation:nthcumulant})
to be different from 0.

Let us summarize our calculation. In order for
(\ref{equation:nonvanishing}) to hold we require: $i$ is
constant on the cycles of $\gamma\pi^{-1}$; each point of
$\{1, k(1) + 1, \dots, k(1) + \cdots + k(n-1) + 1\}$ is in a
separate cycle of $\gamma\pi^{-1}$; and $i$ is constant on
the union of the cycles of $\gamma\pi^{-1}$ containing the
points of $\{1, k(1) + 1, \dots, k(1) + \cdots + k(n-1) +
1\}$. This leaves $\#(\gamma\pi^{-1}) -n+1 $ cycles on which
we can arbitrarily choose values of $i$ (recall that
$\#(\gamma\pi^{-1})$ denotes the number of cycles of
$\gamma\pi^{-1}$). Thus the number of choices for $i$ is
$N^d$ where $d = \#(\gamma\pi^{-1}) -n + 1$. (See Figure 1.)
So if we sum for a fixed such $\pi$ over all free indices
(each choice of them will give the same contribution,
because the cyclic cumulants of the $a_{ij}$ do not depend
on the actual choice of the indices) then we get altogether
for such a $\pi$ the contribution
\begin{align*}
N^{d} \prod_{V\in\pi} \frac {\kk_{\vert V\vert}}{N^{\vert V\vert -1}} & =
N^{d + \#(\pi) - |V_1| - \cdots - |V_{\#(\pi)}|}
\prod_{V\in\pi} {\kk_{\vert V\vert}} \\
& =
N^{d + \#(\pi) -k} \kk_{\pi}(a,\dots,a) \\
& =
N^{-n+2} \kk_{\pi}(a,\dots,a),
\end{align*}
where $d + \#(\pi) -k = \#(\gamma\pi^{-1}) + \#(\pi) -n -k
+1 = -n + 2$ because $\#(\pi) + \#(\gamma\pi^{-1}) = k + 1$
(see \cite[Exercise 9.23]{NS}).

But carrying out the sum over the $\pi$ is now the same as
calculating cumulants of powers of $a$.  So finally we get
the simple result
\begin{equation}\label{eq:cumulants}
\kk_n(\Tr(A^{k(1)}),\dots,\Tr(A^{k(n)}))=N^{-n+2}
\kk_n(a^{k(1)},\dots,a^{k(n)}).
\end{equation}
One should note that, compared to the case of classical
random matrices, there are no subleading orders.  Thus the
limit $N\to\infty$ does not produce any new feature and
contains essentially the same information as the random
variable $a$, i.e., the case $N=1$. In this sense, these
quantum random matrices are less interesting from the point
of view of fluctuations than their classical
counterparts. Still, let us elaborate a bit on what happens
with respect to fluctuations.

First, it is clear from \eqref{eq:cumulants} that all
cumulants of higher order than 2 go to zero, and thus each
centered trace of a power of $A$ goes to a semicircular
element.  The covariance between two such traces of powers
is (actually for any $N$) given by
$$\kk_2(\Tr(A^p),\Tr(A^q))=\kk_2(a^p,a^q).$$
Since those fluctuations depend on the distribution of $a$
we do not have universality for the fluctuations in the
quantum case.

Let us finally also check whether there is some kind of
analogue of the Schwinger-Dyson equations.  We put
$$\ff_1(p(x)):=\lim_{N\to\infty} \kk_1(\tr(p(A)))=\kk_1(p(a))=\ff(p(a))$$
and
$$\ff_2(p(x),q(x)):=\lim_{N\to\infty}\kk_2(\Tr(p(A)),\Tr(q(A)))=\kk_2(p(a),q(a)).$$

Since $\ff_1$ captures just the information about the
distribution of $a$, the first order equation is nothing
else but the definition of the conjugate variable $\xi$ for
$a$, namely for this we just have the equation
$$\ff_1(\xi p(x))=\ff(\xi p(a))=\ff\otimes\ff (\partial p(a))=\ff_1\otimes\ff_1(\partial p(x)).$$
For the second order we have 
$$
\ff_2(\xi p(x), q(x))=\kk_2(\xi p(a), q(a))$$
which yields, by using again the formula for free cumulants
with products as arguments, the following kind of linear
analogue of \eqref{SD:second}:
\begin{equation}
\ff_2(\xi p(x), q(x))=
\ff_2\bigl( \ff\otimes\id (\partial p(x)),q(x)\bigr)+\ff_1\otimes\ff_1\bigl(p(x)\otimes 1\cdot \partial q(x)\bigr).
\end{equation}


\begin{thebibliography}{10}

\bibitem{AJM} J. Ambj\o rn, J. Jurkiewicz, and Y. Makeenko,
  Multiloop correlators for two-dimensional quantum gravity,
  \textit{Phys. Lett.}, \textbf{251B} (1990), 517-524.

\bibitem{E} B. Eynard, \textit{Random Matrices}, Cours de
  Physique Th\'eorique de Saclay, 2000.

\bibitem{CS} S. Curran and R. Speicher, Quantum invariant
  families of matrices in free
  probability, \textit{J. Funct. Anal.}, \textbf{261} (2011), 897-933.

\bibitem{G1} A. Guionnet, \textit{Large Random Matrices:
  Lectures on Macroscopic Asymptotics $($Saint-Flour
  2006$)$}, Springer LNM 1957, Springer, 2009.
  
\bibitem{Gui} A. Guionnet and E. Maurel-Segala, Second order
  asymptotics for matrix models, \textit{Ann. Probab.},
  \textbf{35}, (2007), 2160--2212.
  
\bibitem{J} K. Johansson, On fluctuations of eigenvalues of
  random Hermitian matrices, \textit{Duke Math. J.},
  \textbf{91} (1998), 151-204.  

\bibitem{KMS} T. Kusalik, J. Mingo, and R. Speicher,
  Orthogonal polynomials and fluctuations of random
  matrices, \textit{J. Reine Angew. Math. $($Crelle's
    J.$)$}, \textbf{604} (2007), 1-46.

\bibitem{MSS} J. Mingo, P. Sniady, and R. Speicher, Second
  order freeness and fluctuations of random matrices;
  II. Unitary random matrices, \textit{Adv. Math.},
  \textbf{209} (2007), 212-240.

\bibitem{MST} J. Mingo, R. Speicher, and E. Tan, Second
  Order Cumulants of Products,
  \textit{Trans. Amer. Math. Soc.}, \textbf{361} (2009),
  4751-4781.

\bibitem{NSS} A. Nica, D. Shlyakhtenko, and R. Speicher,
  $R$-cyclic families of random matrices in free
  probability, \textit{J. Funct. Anal.}, \textbf{188}
  (2002), 227-271.

\bibitem{NS} A. Nica and R. Speicher, \textit{Lectures on
  the combinatorics of free probability}, Cambridge
  University Press, 2006.

\bibitem{P} D. Politzer, Random-matrix description of the
  distribution of mesoscopic conductance,
  \textit{Phys. Rev. B}, \textbf{40} (1989), 11917-11919.

\bibitem{V-cyclic} D. Voiculescu, A note on cyclic gradient,
  \textit{Indiana Univ. Math. J.}, \textbf{49} (2000),
  837-841.

\bibitem{VDN} D. Voiculescu, K. Dykema, and A. Nica,
  \textit{Free Probability Theory}, Amer. Math. Soc., 1992.
\end{thebibliography}
\end{document}